\documentclass[12pt]{article}
\usepackage{amssymb, amsmath, amsfonts}
\usepackage{mathrsfs}

\newtheorem{theorem}{Theorem}
\newtheorem{corollary}{Corollary}

\begin{document}

\begin{center}
{\Large {\bf Central limit theorem related to MDR-method}}
%the multifactor dimensionality reduction method}
\end{center}
\begin{center}
\large{\bf Alexander Bulinski}\footnote{Faculty of Mathematics and Mechanics,
Lomonosov Moscow State University, Moscow 119991, Russia.\\
E-mail: \;bulinski@mech.math.msu.su}
$\!\!^,$\footnote{The work is partially supported by RFBR grant 13-01-00612.}
\end{center}
\vskip1cm

In many medical and biological investigations, including
genetics, it is typical to handle high dimensional data which can be
viewed as a set of values of some factors and a binary response
variable. For instance, the response variable can describe the state
of a patient health and one often assumes that it depends only on
some part of factors.  An important problem is to determine
collections of significant factors. In this regard we turn to the
MDR-method introduced by M.Ritchie and coauthors. Our recent paper
provided the necessary and sufficient conditions for strong
consistency of estimates of the prediction error employing the
$K$-fold cross-validation and an arbitrary penalty function. Here we
introduce the regularized versions of the mentioned estimates and
prove for them the multidimensional CLT. Statistical variants of the
CLT involving self-normalization are discussed as well.


\vskip0.3cm

\vskip1cm
\section{Introduction}
\label{sec:1}

High dimensional data arise naturally in a number of experiments.
Very often such data are viewed as the values of some factors
$X_1,\ldots,X_n$ and the corresponding response variable $Y$. For
example, in medical studies such response variable $Y$ can describe
the health state (e.g., $Y=1$ or $Y=-1$  mean ``sick'' or
``healthy'') and  $X_1,\ldots,X_m$ and $X_{m+1},\ldots,X_n$ are
genetic and non-genetic factors, respectively. Usually $X_i$ ($1\leq
i\leq m$) characterizes a single nucleotide polymorphism (SNP), i.e.
a certain change of nucleotide bases adenine, cytosine, thymine and
guanine (these genetic notions can be found, e.g., in \cite{BHFH})
in a specified segment of DNA molecule. In this case one considers
$X_i$ with three values, for instance, $0,1$ and $2$ (see, e.g.,
\cite{BBSSYBST}). It is convenient to suppose that other $X_i$
($m+1\leq i\leq n$)  take values in $\{0,1,2\}$ as well. For
example, the range of blood pressure can be partitioned into zones
of low, normal and high values. However, further we will suppose
that all factors take values in arbitrary finite set. The binary
response variable can also appear in pharmacological experiments
where $Y=1$ means that the medicament is efficient and $Y=-1$
otherwise.

A challenging problem is to find the genetic and non-genetic (or
environmental) factors which could increase the risk of complex
diseases such as diabetes, myocardial infarction and others. Now the
most part of  specialists share the paradigm that in contrast to
simple disease (such as sickle anemia) certain combinations of the
``damages'' of the DNA molecule could be responsible for provoking
the complex disease whereas the single mutations need not have
dangerous effects (see, e.g., \cite{SR}). The important research
domain called the {\it genome-wide association studies} (GWAS)
inspires development of new methods for handling large massives of
biostatistical data. Here we will continue our treatment of the {\it
multifactor dimensionality reduction} (MDR) method introduced by
M.Ritchie et al. \cite{Ritch}. The idea of this method goes back to
the Michalski algorithm. A comprehensive survey concerning the MDR
method is provided in \cite{RM}, on subsequent modifications and
applications see, e.g., \cite{ET}, \cite{Gui}  -- \cite{OLK},
\cite{VWM} and \cite{Winham}. Other complementary methods applied in
GWAS are discussed, e.g., in \cite{BBSSYBST}, there one can find
further references.

In \cite{B} the basis for application of the MDR-method was proposed
when one uses an arbitrary penalty function to describe the
prediction error of the binary response variable by means of a
function in factors. The goal of the present paper is to establish
the new multidimensional central limit theorem (CLT) for statistics
which permit to justify the optimal choice of a subcollection of the
explanatory variables.

\section{Auxiliary results}
\label{sec:2}

Let $X=(X_1,\ldots,X_n)$ be a random vector with components
$X_i:\Omega \to \{0,1,\ldots,q\}$ where $i=1,\ldots, n$ ($q$, $n$ are
positive integers). Thus, $X$ takes values in
$\mathbb{X}=\{0,1,\ldots,q\}^n$. Introduce a random (response)
variable $Y:\Omega \to \{-1,1\}$, non-random function
${f:\mathbb{X}\to \{-1,1\}}$ and a penalty function ${\psi:
\{-1,1\}\to \mathbb{R}_+}$ (the trivial case $\psi \equiv 0$ is
excluded). The quality of approximation of $Y$ by $f(X)$ is defined
as follows
\begin{equation}\label{-1}
Err(f):={\sf E} |Y-f(X)|\psi(Y).
\end{equation}
Set $M=\{x\in\mathbb{X}:{\sf P}(X=x)>0\}$ and $$F(x)= \psi(-1){\sf P}(Y=-1|X=x)-\psi(1){\sf P}(Y=1|X=x),\;\;x\in M.$$

It is not difficult to show (see \cite{B}) that the collection of
{\it optimal functions}, i.e. all functions $f:\mathbb{X}\to
\{-1,1\}$ which are solutions of the problem $Err(f)\to$ {\it inf},
has the form
\begin{equation}\label{0a}
f= \mathbb{I}\{A\}- \mathbb{I}\{\overline{A}\},\;\;A\in \mathcal{A},
\end{equation}
$\mathbb{I}\{A\}$ stands for an indicator of $A$
($\mathbb{I}\{\varnothing\}:=0$) and $\mathcal{A}$ consists of
sets$$A=\{x\in M: F(x)<0\}\cup B\cup C.$$ Here $B$ is an arbitrary
subset of $\{x\in M: F(x)=0\}$ and $C$ is any subset of
$\overline{M}:=\mathbb{X}\setminus M$. If we take $A^*=\{x\in M:
F(x)<0\}$, then $A^*$  has the minimal cardinality among all subsets
of $\mathcal{A}$. In view of the relation $\psi(-1)+\psi(1)\neq 0$
we have
\begin{equation}\label{1}
A^*=\{x\in M: {\sf P}(Y=1|X=x) >
\gamma(\psi)\},\;\;\;\;\gamma(\psi):=\psi(-1)/(\psi(-1)+\psi(1)).
\end{equation}
If $\psi(1)=0$ then $A^*=\varnothing$. If $\psi(1)\neq 0$ and
$\psi(-1)/\psi(1)= a$ where $a\in \mathbb{R}_+$ then
$A^*=\left\{x\in M: {\sf P}(Y=1|X=x) > a/(1+a)\right\}$. Note that
we can rewrite \eqref{-1} as follows
$$
Err(f)= 2\sum_{y\in\{-1,1\}}\psi(y){\sf P}(Y=y,f(X)\neq y).
$$
The value $Err(f)$ is unknown as we do not know the law of a random
vector $(X,Y)$. Thus, statistical inference on the quality of
approximation of $Y$ by means of $f(X)$ is based on the estimate of
$Err(f)$.

Let $\xi^1,\xi^2,\ldots$ be i.i.d. random vectors with the same law
as a vector $(X,Y)$. For $N\in\mathbb{N}$ set
$\xi_N=\{\xi^1,\ldots,\xi^N\}$. To approximate $Err(f)$, as $N\to
\infty$, we will use a
 {\it prediction algorithm}. It involves a function
$f_{PA}=f_{PA}(x,\xi_N)$ with values $\{-1,1\}$ which is defined for
$x\in \mathbb{X}$ and $\xi_N$. In fact we use a {\it family} of
functions $f_{PA}(x,v_m)$ defined for $x\in \mathbb{X}$ and $v_m\in
\mathbb{V}_m$ where $\mathbb{V}_m :=(\mathbb{X} \times \{-1,1\})^m$,
$m\in\mathbb{N}$, $m\leq N$. To simplify the notation we write
$f_{PA}(x,v_m)$ instead of $f_{PA}^m(x,v_m)$. For $S\subset
\{1,\ldots,N\}$ ("$\subset$"\, means non-strict inclusion
"$\subseteq$"\,) put $\xi_N(S)= \{{\xi}^j,j\in S\}$ and
${\overline{S}:=\{1,\ldots,N\}\setminus S}$. For $K\in\mathbb{N}$
($K>1$) introduce a partition of $\{1,\ldots,N\}$ formed by subsets
$$
S_k(N)= \{(k-1)[N/K]+1,\ldots,k[N/K]\mathbb{I}\{k<K\} +
N\mathbb{I}\{k=K\}\}, \;\;k=1,\ldots,K,
$$
here $[b]$ is the integer part of a number $b\in \mathbb{R}$.
Generalizing  \cite{BBSSYBST} we can construct an estimate of
$Err(f)$ using a sample $\xi_N$, a prediction algorithm with
$f_{PA}$ and $K$-fold cross-validation where $K \in\mathbb{N}$,
$K>1$ (on cross-validation see, e.g., \cite{AC}). Namely, let
\begin{equation}\label{5}
\widehat{E}rr_K(f_{PA},\xi_N):= 2\!\!\!\!\!\sum_{y\in\{-1,1\}}\!
\frac{1}{K}\sum_{k=1}^K\sum_{j\in S_k(N)}
\!\!\!\!\!\frac{\widehat{\psi}(y,S_k(N)) \mathbb{I} \{Y^j\!=\!y,
f_{PA}(X^j,\xi_N(\overline{S_k(N)}))\!\neq\! y\}} {\sharp S_k(N)}.
\end{equation}
For each $k=1,\ldots,K$, random variables $\widehat{\psi}(y,S_k(N))$
denote strongly consistent estimates (as $N\to \infty$) of
$\psi(y)$, $y\in \{-1,1\}$, constructed from data $\{Y^j,j\in
S_k(N)\}$, and $\sharp S$ stands for a finite set $S$ cardinality.
We call $\widehat{E}rr_K(f_{PA},\xi_N)$ an {\it estimated prediction
error}.

The following theorem giving a criterion of validity of the relation
\begin{equation}\label{5a}
\widehat{E}rr_K(f_{PA},\xi_N)\to Err(f)\;\;\mbox{a.s.},\;\;N\to
\infty,
\end{equation}
was established in \cite{B} (further on a sum over empty set is
equal to 0 as usual).

\begin{theorem}\label{Th1}
Let $f_{PA}$ define a prediction algorithm for a function
$f\!:\mathbb{X} \to \{-1,1\}$. Assume that there exists such set
$U\subset\mathbb{X}$ that for each $x\in U$ and any $k=1,\ldots,K$
one has
\begin{equation}\label{8}
f_{PA}(x,\xi_N(\overline{S_k(N)}))\to f(x)\;\; \mbox{a.s.},\;\;N\to
\infty.
\end{equation}
Then \eqref{5a} is valid if and only if,  for $N\to \infty$,
\begin{equation}\label{16}
\sum_{k=1}^K\bigl(\!\sum_{x\in\mathbb{X}^+}\!\!
\mathbb{I}\{f_{PA}(x,\xi_N(\overline{S_k(N)}))\!=\!-1\}\!L(x)
\!-\!\sum_{x\in\mathbb{X}^-}\!\!\mathbb{I}\{f_{PA}(x,\xi_N(\overline{S_k(N)}))\!=\!1\}\!L(x)\bigr)
\to 0\; \mbox{a.s.}
\end{equation}
Here $\mathbb{X}^+:=(\mathbb{X}\setminus U) \cap \{x\in M:f(x)=1\}$,
$\mathbb{X}^-:=(\mathbb{X}\setminus U) \cap \{x\in M:f(x)=-1\}$ and
$$
L(x)= \psi(1){\sf P}(X=x,Y=1) - \psi(-1){\sf P}(X=x,Y=-1),\;\;x\in
\mathbb{X}.
$$
\end{theorem}

The sense of this result is the following. It shows that one has to
demand condition \eqref{16} outside the set $U$ (i.e. outside the
set where $f_{PA}$ provides the a.s. approximation of $f$) to obtain
\eqref{5a}.

\begin{corollary}[\rm{\cite{B}}]
Let, for a function $f\!:\mathbb{X} \to \{-1,1\}$, a prediction
algorithm be defined by $f_{PA}$. Suppose that there exists a set
$U\subset\mathbb{X}$ such that for each $x\in U$ and any
$k=1,\ldots,K$ relation \eqref{8} is true. If
$$
L(x)=0\;\;\mbox{for}\;\;x\in (\mathbb{X}\setminus U)\cap M
$$
then \eqref{5a} is satisfied.
\end{corollary}

Note also that Remark 4 from \cite{B} explains why the  choice of a
penalty function  proposed by Velez et al. \cite{VWM}:
\begin{equation}\label{velez}
\psi(y) = c({\sf P}(Y=y))^{-1}, \;\;y\in\{-1,1\},\;\;c>0,
\end{equation}
is natural. Further discussion and examples can be found in
\cite{B}.

\section{Main results and proofs}
\label{sec:3}

In many situations it is reasonable to suppose that the response
variable $Y$ depends only on subcollection $X_{k_1},\ldots,X_{k_r}$
of the explanatory variables,   $\{k_1,\ldots,k_r\}$ being a subset
of $\{1,\ldots,n\}$. It means that for any $x\in M$
\begin{equation}\label{mdr}
{\sf P}(Y=1|X_1=x_1,\ldots,X_n=x_n)= {\sf P}(Y=1|X_{k_1}=x_{k_1},\ldots,X_{k_r}=x_{k_r}).
\end{equation}
In the framework of the complex disease analysis it is natural to
assume that only part of the risk factors could provoke this disease
and the impact of others can be neglected. Any collection
$\{k_1,\ldots,k_r\}$ implying \eqref{mdr} is called {\it
significant}. Evidently if $\{k_1,\ldots,k_r\}$ is signifi\-cant
then any collection $\{m_1,\ldots,m_i\}$ such that
$\{k_1,\ldots,k_r\}\subset \{m_1,\ldots,m_i\}$ is significant as
well. For a set $D\subset \mathbb{X}$ let $\pi_{k_1,\ldots,k_r}D:=
\{u=(x_{k_1},\ldots,x_{k_r}): x=(x_1,\ldots,x_n)\in D\}$. For $B\in
\mathbb{X}_r$ where $\mathbb{X}_r:=\{0,1,\ldots,q\}^r$ define in
$\mathbb{X}=\mathbb{X}_n$ a cylinder
$$
C_{k_1,\ldots,k_r}(B):=\{x=(x_1,\ldots,x_n)\in
\mathbb{X}:(x_{k_1},\ldots,x_{k_r})\in B\}.
$$
For $B=\{u\}$ where $u=(u_1,\ldots,u_r)\in \mathbb{X}_r$ we write $
C_{k_1,\ldots,k_r}(u)$ instead of $C_{k_1,\ldots,k_r}(\{u\})$.
Obviously
$$
{\sf P}(Y=1|X_{k_1}=x_{k_1},\ldots,X_{k_r}=x_{k_r})\equiv {\sf P}(Y=1|X\in C_{k_1,\ldots,k_r}(u)),
$$
here
\begin{equation}\label{ux}
u=\pi_{k_1,\ldots,k_r}\{x\}, \;\;\mbox{i.e.}\;\;u_i=x_{k_i},
\;\;i=1,\ldots,r.
\end{equation}

For $C\subset \mathbb{X}$, $N\in\mathbb{N}$ and $W_N\subset
\{1,\ldots,N\}$ set
\begin{equation}\label{estcp}
\widehat{{\sf P}}_{W_N}(Y=1|X\in C):= \frac{\sum_{j\in
W_N}\mathbb{I}\{Y^j=1,X^j \in C\}}{\sum_{j\in W_N}
\mathbb{I}\{X^j\in C\}}.
\end{equation}
When $C=\mathbb{X}$ we write simply $\widehat{{\sf P}}_{W_N}(Y=1)$
in \eqref{estcp}. According to the {\it strong law of large numbers
for arrays} (SLLNA), see, e.g., \cite{TH}, for any $C\subset
\mathbb{X}$ with ${\sf P}(X\in C)>0$
$$
\widehat{{\sf P}}_{W_N}(Y=1|X\in C)\to {\sf P}(Y=1|X\in
C)\;\;\mbox{a.s.}, \;\;\sharp W_N\to \infty,\;\; N\to \infty.
$$

If \eqref{mdr} is valid then the optimal function $f^*$ defined by
\eqref{0a} with $A=A^*$ introduced in \eqref{1} has the form
\begin{equation}\label{mdr1}
f^{k_1,\ldots,k_r}(x)=
\begin{cases}
\;\;\,1, &\mbox{if}\quad{\sf P}(Y=1|X\in C_{k_1,\ldots,k_r}(u))>\gamma(\psi)\;\;\mbox{and}\;\;x\in M,\\
-1, &\mbox{otherwise},
\end{cases}
\end{equation}
here $u$ and $x$ satisfy \eqref{ux} (${\sf P}(X\in
C_{k_1,\ldots,k_r}(u))\geq {\sf P}(X=x)>0$ as $x\in M$). Hence, for
each significant $\{k_1,\ldots,k_r\}\subset \{1,\ldots,n\}$ and any
$\{m_1,\ldots,m_r\}\subset \{1,\ldots,n\}$ one has
\begin{equation}\label{optim}
 Err(f^{k_1,\ldots,k_r})\leq Err(f^{m_1,\ldots,m_r}).
\end{equation}

For arbitrary $\{m_1,\ldots,m_r\}\subset \{1,\ldots,n\}$, $x\in
\mathbb{X}$,  $u=\pi_{m_1,\ldots,m_r}\{x\}$ and a penalty function
$\psi$ we consider the prediction algorithm with a function
$\widehat{f}_{PA}^{m_1,\ldots,m_r}$ such that
\begin{equation}\label{mdr2}
\widehat{f}_{PA}^{m_1,\ldots,m_r}(x,\xi_N(W_N))=
\begin{cases}
\;\;\,1, &\widehat{{\sf P}}_{W_N}(Y=1|X\in
C_{m_1,\ldots,m_r}(u))>\widehat{\gamma}_{W_N}(\psi),\;\;
x\in M,\\
-1, &\mbox{otherwise},
\end{cases}
\end{equation}
here $\widehat{\gamma}_{W_N}(\psi)$ is a strongly consistent
estimate of $\gamma(\psi)$ constructed by means of $\xi_N(W_N)$.
Introduce
\begin{equation}\label{choice_u}
U:=\{x\in M: {\sf P}(Y=1|X_{m_1}=x_{m_1},\ldots,X_{m_r}
=x_{m_r})\neq\gamma(\psi)\}.
\end{equation}
Using Corollary 1 (and in view of Examples 1 and 2 of \cite{B}) we
conclude that for any $\{m_1,\ldots,m_r\}\subset \{1,\ldots,n\}$
\begin{equation}\label{prin_eq}
\widehat{E}rr_K(\widehat{f}_{PA}^{m_1,\ldots,m_r},\xi_N)\to
Err(f^{m_1,\ldots,m_r}) \;\;\mbox{a.s.},\;\;N\to \infty.
\end{equation}

For each
 $\varepsilon>0$, any significant
collection
 $\{k_1,\ldots,k_r\}\subset \{1,\ldots,n\}$ and arbitrary set
  $\{m_1,\ldots,m_r\}\subset \{1,\ldots,n\}$ due to relations \eqref{optim} and \eqref{prin_eq} one has
\begin{equation}\label{mdr3}
\widehat{E}rr_K(\widehat{f}_{PA}^{k_1,\ldots,k_r},\xi_N)\leq
\widehat{E}rr_K(\widehat{f}_{PA}^{m_1,\ldots,m_r},\xi_N)+\varepsilon\;\;\mbox{a.s.}
\end{equation}
when $N$ is large enough.

Thus, for a given $r=1,\ldots,n-1$, according to \eqref{mdr3} we
come to the following conclusion. It is natural to choose among
factors $X_1,\ldots,X_n$ a collection $X_{k_1},\ldots,X_{k_r}$
leading to the smallest estimated prediction error
${\widehat{E}rr_K(\widehat{f}_{PA}^{k_1,\ldots,k_r},\xi_N)}$. After
that it is desirable to apply the permutation tests (see, e.g.,
\cite{BBSSYBST} and \cite{GLMP}) for validation of the prediction
power of selected factors. We do not tackle here the choice of $r$,
some recommendations can be found in \cite{RM}. Note also in passing
that a nontrivial problem is to estimate the importance
 of various collections of factors, see, e.g.,
\cite{SR}. \vskip0.2cm {\bf Remark 1.} It is essential that for each
$\{m_1,\ldots,m_r\} \subset \{1,\ldots,n\}$ we have strongly
consistent estimates of $Err(f^{m_1,\ldots,m_r})$. So to compare
these estimates we can use the subset of $\Omega$ having probability
one. If we had only the convergence in probability instead of a.s.
convergence in \eqref{prin_eq} then to compare different
${\widehat{E}rr_K(\widehat{f}_{PA}^{m_1,\ldots,m_r},\xi_N)}$ one
should take into account the Bonferroni corrections for all subsets
$\{m_1,\ldots,m_r\}$ of $\{1,\ldots,n\}$.

Further on we consider a function $\psi$ having the form
\eqref{velez}. In view of \eqref{1} w.l.g. we can assume that $c=1$
in \eqref{velez}. In this case $\gamma(\psi)= {\sf P}(Y=1).$
Introduce events
$$A_{N,k}(y)=\{Y^j=-y,\;\; j\in S_k(N)\},\;\;N\in \mathbb{N},\;\; k=1,\ldots,K,\;\;y\in\{-1,1\},
$$
and random variables
$$
\widehat{\psi}_{N,k}(y):= \frac{\mathbb{I}\{\overline{A_{N,k}(y)}\}}
{\widehat{\sf P}_{S_k(N)}(Y=y)},
$$
where we write $\widehat{\psi}_{N,k}(y)$ instead of $\widehat{\psi}(y,S_k(N))$.
Trivial cases ${\sf P}(Y=y)\in \{0,1\}$ are excluded and we
formally set $0/0:=0$.  Then
\begin{equation}\label{imp_converg}
\widehat{\psi}_{N,k}(y)- \psi(y) = \frac{{\sf P}(Y=y)- \widehat{\sf P}_{S_k(N)}(Y=y)}{\widehat{\sf P}_{S_k(N)}(Y=y){\sf P}(Y=y)}
\mathbb{I}\{\overline{A_{N,k}(y)}\} -\frac{1}{{\sf P}(Y=y)}\mathbb{I}\{A_{N,k}(y)\}.
\end{equation}
Clearly,
\begin{equation}\label{imp_conv_b}
\mathbb{I}\{A_{N,k}(y)\}\to
0\;\;\mbox{a.s.},\;\;N\to \infty,
\end{equation}
and the following relation is true
\begin{equation}\label{imp_conv_a}
\frac{\mathbb{I}\{\overline{A_{N,k}(y)}\}}{\widehat{\sf P}_{S_k(N)}(Y=y)}
 \to \frac{1}{{\sf P}(Y=y)}\;\;\mbox{a.s.},\;\;N\to \infty.
\end{equation}
Therefore, by virtue of \eqref{imp_converg} -- \eqref{imp_conv_a} we
have that for $y\in\{-1,1\}$ and $k=1,\ldots,K$
\begin{equation}\label{eq_fin2}
\widehat{\psi}_{N,k}(y) - \psi(y)\to
0\;\;\mbox{a.s.},\;\;N\to \infty.
\end{equation}

Let $\{m_1,\ldots,m_r\}\subset \{1,\ldots,n\}$. We define the
functions which can be viewed as the {\it regularized versions} of
the estimates $\widehat{f}_{PA}^{m_1,\ldots,m_r}$ of
$f^{m_1,\ldots,m_r}$ (see \eqref{mdr2} and \eqref{mdr1}).
Namely, for $W_N\subset \{1,\ldots,N\}$, $N\in\mathbb{N}$, and
$\varepsilon=(\varepsilon_N)_{N\in\mathbb{N}}$ where non-random
positive $\varepsilon_N\to 0$, as $N\to \infty$, put
$$\widehat{f}_{PA,\varepsilon}^{m_1,\ldots,m_r}(x,\xi_N(W_N))=
\begin{cases}
\;\,1, \!\!\!\!&\widehat{\sf P}_{W_N}(Y\!=\!1|X\in C_{m_1,\ldots,m_r}(u))>\widehat\gamma_{W_N}(\psi)+\varepsilon_N, \;x\in M,\\
-1, \!\!&\mbox{otherwise},
\end{cases}
$$
where $u=\pi_{m_1,\ldots,m_r}\{x\}$. Regularization of
$\widehat{f}^{m_1,\ldots,m_r}_{PA}$ means that instead of  the
thre\-shold $\widehat{\gamma}_{W_N}(\psi)$ we use
$\widehat{\gamma}_{W_N}(\psi)+\varepsilon_N$.

Take now $U$ appearing in \eqref{choice_u}. Applying Corollary 1
once again (and in view of Examples 1 and 2 of \cite{B}) we can
claim that the statements which are analogous to \eqref{prin_eq} and
\eqref{mdr3} are valid for the regularized versions of the estimates
introduced above. Now we turn to the principle results, namely,
central limit theorems.

\begin{theorem}\label{Th2}
Let $\varepsilon_N\to 0$ and $N^{1/2}\varepsilon_N \to \infty$  as
$N\to \infty$. Then, for each $K\in \mathbb{N}$, any subset
$\{m_1,\ldots m_r\}$ of $\{1,\ldots,n\}$, the corresponding function
$f=f^{m_1,\ldots,m_r}$ and prediction algorithm defined by
$f_{PA}=\widehat{f}_{PA,\varepsilon}^{m_1,\ldots,m_r}$, the
following relation holds:
\begin{equation}\label{CLT}
\sqrt{N}(\widehat{E}rr_K(f_{PA},\xi_N) - Err(f))
\stackrel{law}\longrightarrow Z\sim N(0,\sigma^2),\;\;N\to \infty,
\end{equation}
where $\sigma^2$ is variance of the random variable
\begin{equation}\label{abc}
V= 2\sum_{y\in \{-1,1\}} \frac{\mathbb{I}\{Y=y\}}{{\sf P}(Y=y)}
\left(\mathbb{I}\{f(X)\neq y\} - {\sf P}(f(X)\neq y|Y=y)\right).
\end{equation}
\end{theorem}

{\bf Proof.}
For a fixed $K\in\mathbb{N}$ and any $N\in\mathbb{N}$ set
$$
T_{N}(f):= \frac{2}{K} \!\sum_{k=1}^K\frac{1}{\sharp
S_k(N)}\sum_{y\in\{-1,1\}}\psi(y) \sum_{j\in S_k(N)}\mathbb{I}
\{Y^j=y, f(X^j)\neq y\},
$$
$$
\widehat{T}_{N}(f):= \frac{2}{K} \!\sum_{k=1}^K\frac{1}{\sharp
S_k(N)}\sum_{y\in\{-1,1\}}\widehat{\psi}_{N,k}(y) \sum_{j\in
S_k(N)}\mathbb{I} \{Y^j=y, f(X^j)\neq y\}.
$$
One has
$$
\widehat{E}rr_K(f_{PA},\xi_N)- Err(f) =
(\widehat{E}rr_K(f_{PA},\xi_N) - \widehat{T}_N(f))
$$
\begin{equation}\label{represent}
+(\widehat{T}_N(f) - T_N(f)) + (T_N(f)-Err(f)).
\end{equation}
First of all we show that
\begin{equation}\label{conv}
\sqrt{N}(\widehat{E}rr_K(f_{PA},\xi_N) -
\widehat{T}_N(f))\stackrel{\sf P}\longrightarrow 0,\;\;N\to \infty.
\end{equation}
For $x\in \mathbb{X}$, $y\in\{-1,1\}$, $k=1,\ldots,K$ and
$N\in\mathbb{N}$ introduce
$$
F_{N,k}(x,y):= \mathbb{I}\{f_{PA}(x,\xi_N(\overline{S_k(N)}))\neq
y\} - \mathbb{I}\{f(x)\neq y\}.
$$
Then
\begin{equation}\label{aa}
\widehat{E}rr_K(f_{PA},\xi_N) - \widehat{T}_{N}(f) = \frac{2}{K}
\!\sum_{k=1}^K \frac{1}{\sharp S_k(N)}
\sum_{y\in\{-1,1\}}\!\!\!\widehat{\psi}_{N,k}(y) \sum_{j\in
S_k(N)}\mathbb{I} \{Y^j=y\}F_{N,k}(X^j,y).
\end{equation}
We define the random variables
$$
B_{N,k}(y):= \frac{1}{\sqrt{\sharp S_k(N)}}\sum_{j\in
S_k(N)}\mathbb{I} \{Y^j=y\}F_{N,k}(X^j,y)
$$
and verify that for each $k=1,\ldots,K$
\begin{equation}\label{conv_pr}
\sum_{y\in\{-1,1\}}\widehat{\psi}_{N,k}(y)B_{N,k}(y) \stackrel{\sf P}\longrightarrow 0,\;\;N\to \infty.
\end{equation}
Clearly \eqref{conv_pr} implies \eqref{conv} in view of \eqref{aa}
as $\sharp S_k(N)= [N/K]$ for $k=1,\ldots,K-1$ and $[N/K]\leq \sharp
S_K(N)<[N/K]+K$. Write $B_{N,k}(y)=
B_{N,k}^{(1)}(y)+B_{N,k}^{(2)}(y)$ where
$$
B_{N,k}^{(1)}(y)\!=\! \frac{1}{\sqrt{\sharp S_k(N)}}\!\!\sum_{j\in
S_k(N)}\!\!\!\!\mathbb{I}\{X^j\!\in\!
U)\}\mathbb{I}\{Y^j=y\}F_{N,k}(X^j,y),
$$
$$
B_{N,k}^{(2)}(y)\!=\! \frac{1}{\sqrt{\sharp S_k(N)}}\!\!\sum_{j\in
S_k(N)}\!\!\!\!\mathbb{I}\{X^j\!\notin \!
U\}\mathbb{I}\{Y^j=y\}F_{N,k}(X^j,y).
$$
Obviously
$$
|B_{N,k}^{(1)}(y)|\leq \sum_{x\in U}\frac{1}{\sqrt{\sharp
S_k(N)}}\sum_{j\in
S_k(N)}|\mathbb{I}\{f_{PA}(x,\xi_N(\overline{S_k(N)}))\neq y\} -
\mathbb{I}\{f(x)\neq y\}|.
$$
Functions $f_{PA}$ and $f$ take values in the set $\{-1,1\}$. Thus,
for any $x\in U$ (where $U$ is defined in \eqref{choice_u}),
$k=1,\ldots,K$ and almost all $\omega \in \Omega$ relation \eqref{8}
ensures the existence of an integer $N_{0}(x,k,\omega)$ such that
$f_{PA}(x,\xi_N(\overline{S_k(N)}))=f(x)$ for $N\geq
N_0(x,k,\omega)$. Hence $B_{N,k}^{(1)}(y)=0$ for any $y$ belonging
to $\{-1,1\}$, each $k=1,\ldots,K$ and almost all $\omega \in
\Omega$ when $N\geq N_{0,k}(\omega)=\max_{x\in U}N_0(x,k,\omega)$.
Evidently, $N_{0,k}<\infty$ a.s., because $\sharp U < \infty$. We
obtain that
\begin{equation}\label{conv_pr_1}
\sum_{y\in\{-1,1\}}\widehat{\psi}_{N,k}(y)B_{N,k}^{(1)}(y) \to
0\;\;\mbox{a.s.},\;\;\;N\to \infty.
\end{equation}
If $U=\mathbb{X}$ then $B_{N,k}^{(2)}(y)=0$ for all $N,k$ and $y$
under consideration. Consequently, \eqref{conv_pr} is valid and
thus, for $U=\mathbb{X}$, relation \eqref{conv} holds. Let now
$U\neq \mathbb{X}$. Then for  $k=1,\ldots,K$ and $N\in\mathbb{N}$
one has
$$
\sum_{y\in\{-1,1\}}\widehat{\psi}_{N,k}(y) B_{N,k}^{(2)}(y)=
\sum_{x\in\mathbb{X}_+}\sum_{y\in\{-1,1\}}H_{N,k}(x,y) +
\sum_{x\in\mathbb{X}_-}\sum_{y\in\{-1,1\}}H_{N,k}(x,y),
$$
here $\mathbb{X}_+=(\mathbb{X}\setminus U)\cap\{x\in
\mathbb{X}:f(x)=1\}$, $\mathbb{X}_-=(\mathbb{X}\setminus
U)\cap\{x\in \mathbb{X}:f(x)=-1\}$ and
$$
H_{N,k}(x,y)\!:=\!\frac{\widehat{\psi}_{N,k}(y)}{\sqrt{\sharp
S_k(N)}}\!\!\!\sum_{j\in S_k(N)}\!\!\!\!\mathbb{I}\{A^j(x,y)\}
(\mathbb{I}\{f_{PA}(x,\xi_N(\overline{S_k(N)}))\!\neq\! y\} -
\mathbb{I}\{f(x)\!\neq \!y\})
$$
where $A^j(x,y)=\{X^j=x,Y^j=y\}$. The definition of $U$ yields that
$\mathbb{X}_+=\varnothing$ and
$$\mathbb{X}_-=\overline{M}\cup\{x\in M:
{\sf P}(Y=1|X_{m_1}=x_{m_1},\ldots,X_{m_r}= x_{m_r}) =
\gamma(\psi)\}.$$ Set
$$\widehat{R}^j_{N,k}(x)=\mathbb{I}\{X^j=x\}(\widehat{\psi}_{N,k}(1)\mathbb{I}\{Y^j=1\}-
\widehat{\psi}_{N,k}(-1)\mathbb{I}\{Y^j=-1\}).$$ It is easily seen
that
$$
\sum_{x\in\mathbb{X}_-}\sum_{y\in\{-1,1\}}H_{N,k}(x,y)=
-\sum_{x\in\mathbb{X}_-}\mathbb{I}\{f_{PA}(x,\xi_N(\overline{S_k(N)}))=1)\}
\!\!\!\sum_{j\in S_k(N)}\frac{\widehat{R}^j_{N,k}(x)}{\sqrt{\sharp
S_k(N)}}.
$$
Note that $\widehat{R}^j_{N,k}(x)=0$ a.s. for all
$x\in\overline{M}$, $k=1,\ldots,K$, $j=1,\ldots,N$ and
$N\in\mathbb{N}$. Let us prove that, for any $x\in
M\cap\mathbb{X}_-$ and $k=1,\ldots,K$,
\begin{equation}\label{conv_probab}
\mathbb{I}\{f_{PA}(x,\xi_N(\overline{S_k(N)}))=1\} \stackrel{\sf P}\longrightarrow 0,\;\;N\to \infty.
\end{equation}
For any $\nu >0$ and $x\in M\cap\mathbb{X}_-$ we have
$$
{\sf P}(\mathbb{I}\{f_{PA}(x,\xi_N(\overline{S_k(N)}))=1\}>\nu)
$$
$$
={\sf P}\left(\widehat{{\sf P}}_{\;\overline{S_k(N)}}(Y=1|X_{m_1}=x_{m_1},\ldots,X_{m_r}=
x_{m_r})>\widehat{\gamma}_{\;\overline{S_k(N)}}(\psi)+
\varepsilon_N\right).
$$
Now we show that, for $k=1,\ldots,K$, this probability tends to $0$
as $N\to \infty$. For ${W_N\subset\{1,\ldots,N\}}$ and $x\in
M\cap\mathbb{X}_-$, put
$$
\Delta_N(W_N,x):={\sf P}\left(\frac{\frac{1}{\sharp W_N}\sum_{j\in
W_N} \eta^j} {\frac{1}{\sharp W_N}\sum_{j\in
W_N}\zeta^j}>\widehat{\gamma}_{W_N}(\psi)+ \varepsilon_N \right)
$$
where $\eta^j=\mathbb{I}\{Y^j
=1,X^j_{m_1}=x_{m_1},\ldots,X^j_{m_r}=x_{m_r}\}$, $\zeta^j=
\mathbb{I}\{X^j_{m_1}=x_{m_1},\ldots,X^j_{m_r}=x_{m_r}\}$,
$j=1,\ldots,N$. Set $p= {\sf P}(X_{m_1}=x_{m_1},\ldots,X_{m_r}=x_{m_r})$. It follows that, for
any $\alpha_N>0$,
\begin{eqnarray}
\!\!& &\quad\quad\quad\quad\quad\quad\quad\quad\quad\quad\quad\quad\quad\quad\quad\Delta_N(W_N,x)\nonumber\\
\!\!&\leq&{\sf P}\Bigl(\frac{\sum_{j\in W_N} \eta^j} {\sum_{j\in
W_N}\zeta^j}>\widehat{\gamma}_{W_N}(\psi)+ \varepsilon_N,
\Bigl|\frac{1}{\sharp W_N} \sum_{j\in W_N}\zeta^j - p\Bigr|<
\alpha_N, \Bigl|\widehat{\gamma}_{W_N}(\psi)-\gamma(\psi)\Bigr|<
\alpha_N\Bigr)\nonumber\\
\!\!&+&{\sf P}\Bigl(\Bigl|\frac{1}{\sharp W_N}\sum_{j\in W_N}
\zeta^j - p\Bigr|\geq \alpha_N\Bigr) +{\sf P}\Bigl(\Bigl|\frac{1}{\sharp W_N}\sum_{j\in W_N}\mathbb{I}\{Y^j=1\}
- {\sf P}(Y=1)\Bigr|\geq \alpha_N\Bigr).\quad\quad\label{ineq_10}
\end{eqnarray}
Due to the Hoeffding inequality
$$
{\sf P}\Bigl(\Bigl|\frac{1}{\sharp W_N}\sum_{j\in W_N} \zeta^j -
p\Bigr|\geq \alpha_N\Bigr) \leq 2 \exp\left\{-2\sharp W_N
\alpha_N^2\right\}=:\delta_N(W_N,\alpha_N).
$$
We have an analogous estimate for the last summand in
\eqref{ineq_10}. Consequently, taking into account that $p>0$ we see
that for all $N$ large enough
$$
\Delta_N(W_N,x) \leq {\sf P}\Bigl(\frac{1}{\sharp W_N}\sum_{j\in
W_N} \eta^j
> (p-\alpha_N)(\gamma(\psi)-\alpha_N+\varepsilon_N)\Bigr)+2\delta_N(W_N,\alpha_N).
$$

Whenever $x\in M\cap\mathbb{X}_-$ one has
$$
{\sf P}(Y=1,X_{m_1}=x_{m_1},\ldots,X_{m_r}=x_{m_r}) = {\sf P}(Y=1){\sf P}(X_{m_1}=x_{m_1},\ldots,X_{m_r}=x_{m_r}),
$$
therefore
$$
\Delta_N(W_N,x)\leq {\sf P}\Bigl(\sum_{j\in W_N} \frac{\eta^j-{\sf E}\eta^j}{\sqrt{\sharp W_N}}
>\sqrt{\sharp W_N} \bigl(p\varepsilon_N-
\alpha_N(\gamma(\psi)+p-\alpha_N+\varepsilon_N)\bigr)\Bigr)
+2\delta_N(W_N,\alpha_N).
$$

The CLT holds for an array ${\{\eta^j, j\in W_N, N\in\mathbb{N}\}}$
consisting of i.i.d. random variables, thus
$$
\frac{1}{\sqrt{\sharp W_N}}\sum_{j\in W_N} (\eta^j-{\sf E}\eta^j)\stackrel{law}\longrightarrow Z\sim N(0,\sigma_0^2),
$$
here $\sigma_0^2= var \mathbb{I}\{Y
=1,X_{m_1}=x_{m_1},\ldots,X_{m_r}=x_{m_r}\}.$ Hence
{$\Delta_N(W_N,x)\to 0$} if, for some $\alpha_N>0$,
\begin{equation}\label{cond_mu}
\alpha_N\sqrt{\sharp W_N}\to \infty,\;\; \varepsilon_N \sqrt{\sharp
W_N}\to \infty,\;\;\alpha_N/\varepsilon_N\to 0\;\;\mbox{as}
\;\;N\to\infty.
\end{equation}
Take $W_N=\overline{S_k(N)}$ with $k=1,\ldots,K$. Then $\sharp
\overline{S_k(N)} \geq (K-1)[N/K]$ for $k=1,\ldots,K$ and we
conclude that \eqref{cond_mu} is satisfied when $\varepsilon_N
N^{1/2}\to \infty$ as $N\to \infty$ if we choose a sequence
$(\alpha_N)_{N\in\mathbb{N}}$ in appropriate way. So, relation
\eqref{conv_probab} is established.

Let
$$R_j(x)=
\mathbb{I}\{X^j=x\}(\psi(1)\mathbb{I}\{Y^j=1\} -
\psi(-1)\mathbb{I}\{Y^j=-1\}),\;\;x\in\mathbb{X},\;\;j\in\mathbb{N}.$$
For all $x\in M\cap\mathbb{X}_-$ one has
$$
\frac{1}{\sqrt{\sharp S_k(N)}}\!\!\!\sum_{j\in
S_k(N)}\widehat{R}^j_{N,k}(x) = \frac{1}{\sqrt{\sharp
S_k(N)}}\!\!\!\sum_{j\in S_k(N)}R_j(x)\;+
$$
$$
+\! \sum_{j\in S_k(N)}\!\! \mathbb{I}\{X^j=x\}
\frac{(\widehat{\psi}_{N,k}(1) - \psi(1))\mathbb{I}\{Y^j=1\} -
(\widehat{\psi}_{N,k}(-1) -
\psi(-1))\mathbb{I}\{Y^j=-1\}}{\sqrt{\sharp S_k(N)}}.
$$
Note that ${\sf E}R_j(x)=0$ for all $j \in \mathbb{N}$ and
$x\in\mathbb{X}_-$. The CLT for an array of i.i.d. random variables
$\{R_j(x), j\in S_k(N), N\in\mathbb{N}\}$ provides that
$$
\frac{1}{\sqrt{\sharp S_k(N)}}\!\!\! \sum_{j\in S_k(N)}R_j(x)
\stackrel{law}\longrightarrow Z_1\sim
N(0,\sigma_1^2(x)),\;\;N\to\infty,
$$
where $\sigma_1^2(x)= var
(\mathbb{I}\{X=x\}(\psi(1)\mathbb{I}\{Y=1\}-\psi(-1)\mathbb{I}\{Y=-1\}))$,
$x\in\mathbb{X}_-$. For each $y\in \{-1,1\}$,
$$
(\widehat{\psi}_{N,k}(y) - \psi(y))\frac{1}{\sqrt{\sharp
S_k(N)}}\!\!\!\sum_{j\in S_k(N)}
\mathbb{I}\{X^j=x\}\mathbb{I}\{Y^j=y\}
$$
$$
=(\widehat{\psi}_{N,k}(y) - \psi(y))\frac{1}{\sqrt{\sharp
S_k(N)}}\!\!\!\sum_{j\in S_k(N)}
(\mathbb{I}\{X^j=x\}\mathbb{I}\{Y^j=y\}-{\sf E}\mathbb{I}\{X^j=x\}\mathbb{I}\{Y^j=y\})
$$
$$
+(\widehat{\psi}_{N,k}(y) - \psi(y))\sqrt{\sharp S_k(N)}{\sf P}(X=x,Y=y).
$$
Due to the CLT
$$
\sum_{j\in S_k(N)}
\!\!\frac{\mathbb{I}\{X^j=x\}\mathbb{I}\{Y^j=y\}-{\sf E}\mathbb{I}\{X^j=x\}\mathbb{I}\{Y^j=y\}}{\sqrt{\sharp S_k(N)}}
\stackrel{law}\longrightarrow Z_2\!\sim\! N(0,\sigma_2^2(x,y))
$$
as $N\to \infty$, where $\sigma_2^2(x,y)= var
\mathbb{I}\{X^j=x,Y^j=y\}$. In view of \eqref{eq_fin2} we have
$$
\frac{\widehat{\psi}_{N,k}(y) - \psi(y)}{\sqrt{\sharp
S_k(N)}}\!\!\!\sum_{j\in S_k(N)}
(\mathbb{I}\{X^j=x\}\mathbb{I}\{Y^j=y\}-{\sf E}\mathbb{I}\{X^j=x\}\mathbb{I}\{Y^j=y\}) \stackrel{\sf P}\longrightarrow 0
$$
as $N\to \infty$. Now we apply \eqref{imp_converg} --
\eqref{imp_conv_a} once again to conclude that
$$
(\widehat{\psi}_{N,k}(y) - \psi(y))\sqrt{\sharp
S_k(N)}\stackrel{law}\longrightarrow Z_3\sim
N(0,\sigma_3^2(y)),\;\;N\to\infty,
$$
with $\sigma_3^2(y)= {\sf P}(Y=-y)({\sf P}(Y=y))^{-3}$. Thus,
\begin{equation}\label{conv_pr_2}
\sum_{y\in\{-1,1\}}\widehat{\psi}_{N,k}(y)
B_{N,k}^{(2)}(y)\stackrel{\sf P}\longrightarrow 0, \;\;N\to\infty.
\end{equation}
Taking into account \eqref{conv_pr_1} and \eqref{conv_pr_2} we come
to \eqref{conv_pr} and consequently to \eqref{conv}.

Now we turn to the study of $\widehat{T}_N(f) - T_N(f)$ appearing in
\eqref{represent}. One has
\begin{eqnarray*}
&
&\quad\quad\quad\quad\quad\quad\quad\quad\quad\quad\sqrt{N}(\widehat{T}_N(f)
- T_N(f))\\
&=&\frac{2\sqrt{N}}{K} \!\sum_{k=1}^K\frac{1}{\sharp
S_k(N)}\sum_{y\in\{-1,1\}}(\widehat{\psi}_{N,k}(y)-\psi(y))
\sum_{j\in S_k(N)}\mathbb{I} \{Y^j=y, f(X^j)\neq y\}.
\end{eqnarray*}
Put $Z^j=\mathbb{I} \{Y^j=y, f(X^j)\neq y\}$, $j=1,\ldots,N$. For
each $k=1,\ldots,K$
$$
\sum_{y\in\{-1,1\}}(\widehat{\psi}_{N,k}(y)
 -\psi(y))
\frac{1}{\sqrt{\sharp S_k(N)}}\sum_{j\in S_k(N)} \mathbb{I} \{Y^j=y,
f(X^j)\neq y\}=
$$
$$
= \sum_{y\in\{-1,1\}}(\widehat{\psi}_{N,k}(y)
 -\psi(y))
\frac{1}{\sqrt{\sharp S_k(N)}}\sum_{j\in S_k(N)}(Z^j-{\sf E}Z^j)
$$
$$
+ \sqrt{\sharp S_k(N)}\sum_{y\in\{-1,1\}}(\widehat{\psi}_{N,k}(y)
 -\psi(y)){\sf P}(Y=y,f(X)\neq
y).
$$
Due to \eqref{eq_fin2} and CLT for an array of $\{Z^j,j\in S_k(N),
N\in\mathbb{N}\}$ we have
$$
\sum_{y\in\{-1,1\}}\!\!\!(\widehat{\psi}_{N,k}(y)
 -\psi(y))
\frac{1}{\sqrt{\sharp S_k(N)}}\!\!\sum_{j\in S_k(N)}(Z^j-{\sf E}Z^j)\stackrel{\sf P}\longrightarrow 0
$$
as $N\to \infty$. Consequently the limit distribution of
$$
\sqrt{N}[(\widehat{T}_N(f) - T_N(f)) + (T_N(f)-Err(f))]
$$
will be the same as for random variables
\begin{equation}\label{imp_equ}
\sqrt{N}[(T_N(f)-Err(f)) + \frac{2}{K}\sum_{k=1}^K
\sum_{y\in\{-1,1\}} (\widehat{\psi}_{N,k}(y)
 -\psi(y)){\sf P}(Y=y,f(X)\neq y)].
\end{equation}
Note that for each $y\in\{-1,1\}$ and $k=1,\ldots,K$
$$
\widehat{\sf P}_{S_k(N)}(Y=y)- {\sf P}(Y=y)\stackrel{\sf P}\longrightarrow 0,
$$
$$\sqrt{\sharp S_k(N)}(\widehat{\sf P}_{S_k(N)}(Y=y)- {\sf P}(Y=y))
\stackrel{law}\longrightarrow Z_4\sim N(0,\sigma^2_4),
$$
as $N\to \infty$, where $\sigma^2_4={\sf P}(Y=-1){\sf P}(Y=1)$.

Now the Slutsky lemma shows that the limit behavior of the random
variables introduced in \eqref{imp_equ} will be the same as for
random variables
$$
\sqrt{N}(T_N(f)-Err(f))
$$
$$- \frac{2\sqrt{N}}{K}\sum_{k=1}^K
\sum_{y\in\{-1,1\}} \frac{(\widehat{\sf P}_{S_k(N)}(Y=y)- {\sf P}(Y=y)){\sf P}(Y=y,f(X)\neq y)}{{\sf P}(Y=y)^2}
$$
$$
= \frac{2\sqrt{N}}{K}\sum_{k=1}^K \sum_{y\in\{-1,1\}}
\frac{1}{\sharp S_k(N)}\sum_{j\in S_k(N)}
\Bigl(\frac{\mathbb{I}\{Y^j=y,f(X^j)\neq y\} - {\sf P}(Y=y,f(X)\neq
y)}{{\sf P}(Y=y)}
$$
$$
- \frac{\mathbb{I}\{Y^j=y\} - {\sf P}(Y=y){\sf P}(Y=y,f(X)\neq
y)}{{\sf P}(Y=y)^2}\Bigr)
$$
$$
=\frac{\sqrt{N}}{K}\sum_{k=1}^K \frac{1}{\sharp S_k(N)}\sum_{j\in
S_k(N)}(V^j-{\sf E}V^j)
$$
where
$$
V^j= \sum_{y\in\{-1,1\}}\frac{2\mathbb{I}\{Y^j=y\}}{{\sf P}(Y=y)}\left(\mathbb{I}\{f(X^j)\neq y\}- \frac{{\sf P}(Y=y,f(X)\neq
y)}{{\sf P}(Y=y)}\right).
$$
For each $k=1,\ldots,K$, the CLT for an array $\{V^j, j\in S_k(N),
N\in\mathbb{N}\}$ of  i.i.d. random variables
 yields the relation
$$
Z_{N,k}:=\frac{1}{\sqrt{\sharp S_k(N)}}\sum_{j\in S_k(N)}(V^j-{\sf
E}V^j)\stackrel{law}\longrightarrow Z \sim N(0,\sigma^2),\;\; N\to
\infty,
$$
where $\sigma^2 = var\, V$ and $V$ was introduced in \eqref{abc}.
Since $Z_{N,1},\ldots,Z_{N,K}$ are independent and
$\sqrt{N}/{\sqrt{\sharp S_k(N)}}\to \sqrt{K}$ for $k=1,\ldots,K$, as
$N\to \infty$, we come to \eqref{CLT}. The proof is complete.
$\square$

Recall that for a sequence of random variables
$(\eta_N)_{N\in\mathbb{N}}$ and a  sequence of positive numbers
$(a_N)_{N\in\mathbb{N}}$  one writes $\eta_N= o_P(a_N)$ if
$\eta_N/a_N\stackrel{\sf P}\longrightarrow 0$, $N\to \infty$.
\vskip0.2cm {\bf Remark 2.} As usual one can view the CLT as a
result describing the exact rate of approximation for random
variables under consideration. Theorem 2 implies that
\begin{equation}\label{rate}
\widehat{E}rr_K(f_{PA},\xi_N) - Err(f) = o_P(a_N),\;\;N\to \infty,
\end{equation}
where $a_N=o(N^{-1/2})$. The last relation is optimal in a sense
whenever $\sigma^{2}>0$, i.e. one cannot take $a_N=O(N^{-1/2})$ in
\eqref{rate}. \vskip0.2cm {\bf Remark 3.} In view of \eqref{estcp}
it is not difficult  to construct the consistent estimates
$\widehat{\sigma}_N$ of unknown $\sigma$ appearing in \eqref{CLT}.
Therefore (if $\sigma^2\neq 0$) we can claim that under conditions
of Theorem 1
$$
\frac{\sqrt{N}}{\widehat{\sigma}_N}(\widehat{E}rr_K(f_{PA},\xi_N) -
Err(f)) \stackrel{law}\longrightarrow \frac{Z}{\sigma}\sim
N(0,1),\;\;N\to \infty.
$$
\vskip0.2cm Now we consider the multidimensional version of Theorem
2. To simplify notation set $\alpha = (m_1,\ldots,m_r)$. We write
$\widehat{f}^{\alpha}_{PA,\varepsilon}$ and $f^{\alpha}$ instead of
$\widehat{f}^{m_1,\ldots,m_r}_{PA,\varepsilon}$ and
$f^{m_1,\ldots,m_r}$, respectively. Employing the Cram\'er--Wold
device and the proof of Theorem 2 we come to the following statement
(as usual we use the column vectors and write $\top$ for
transposition).

\begin{theorem}\label{Th3}
Let $\varepsilon_N\to 0$ and $N^{1/2}\varepsilon_N \to \infty$  as
$N\to \infty$. Then, for each $K\in \mathbb{N}$, any
$\alpha{(i)}=\{m_1^{(i)},\ldots,m_r^{(i)}\}\subset \{1,\ldots,n\}$
where $i=1,\ldots,s$, one has
$$
\sqrt{N}(Z^{(1)}_N,\ldots,Z^{(s)}_N)^{\top}\stackrel{law}\longrightarrow
\mathsf{Z}\sim N(0,C),\;\; N\to \infty.
$$
Here $Z^{(i)}_N=
\widehat{E}rr_K(\widehat{f}^{\alpha{(i)}}_{PA,\varepsilon},\xi_N)
 - Err(f^{\alpha{(i)}})$, $i=1,\ldots,s$, and the
elements of covariance matrix $C=(c_{i,j})$ have the form
$$
c_{i,j}= cov(V(\alpha{(i)}),V(\alpha{(j)})),\;\;i,j=1,\ldots,s,
$$
 the random variables $V(\alpha{(i)})$ being defined in the same way
as $V$ in \eqref{abc} with  $f^{m_1,\ldots,m_r}$ replaced by
$f^{\alpha{(i)}}$.
\end{theorem}

To conclude we note  (see also Remark 3) that one can construct the
consistent estimates $\widehat{C}_N$ of the unknown (nondegenerate)
covariance matrix $C$ to obtain the statistical version of the last
theorem. Namely, under conditions of Theorem 3 the following
relation is valid
$$(\widehat{C}_N)^{-1/2}(Z^{(1)}_N,\ldots,Z^{(s)}_N)^{\top}\stackrel{law}\longrightarrow
C^{-1/2}\mathsf{Z}\sim N(0,I),\;\; N\to \infty,
$$
where $I$ stands for the unit matrix of order $s$.

\end{document}